\documentclass[pdflatex,sn-mathphys]{sn-jnl}

\jyear{2021}%

\theoremstyle{thmstyleone}%
\theoremstyle{thmstyletwo}%

\theoremstyle{thmstylethree}%

\begin{document}

\title[Mean field game for modeling of COVID-19 spread]{Mean field game for modeling of COVID-19 spread}


\author*[1]{\fnm{Viktoriya} \sur{Petrakova}}\email{vika-svetlakova@yandex.ru}

\author[2,3]{\fnm{Olga} \sur{Krivorotko}}\email{krivorotko.olya@mail.ru}



\affil*[1]{ \orgname{Institute of computational Modeling SB RAS}, \orgaddress{\street{50/44 Akademgorodok str.}, \city{Krasnoyarsk}, \postcode{660036}, \country{Russia}}}

\affil[2]{\orgname{Institute of Computational Mathematics and Mathematical Geophysics of SB RAS}, \orgaddress{\street{6 Ac. Lavrentieva ave.}, \city{Novosibirsk}, \postcode{630090}, \country{Russia}}}

\affil[3]{\orgname{Novosibirsk State University}, \orgaddress{\street{2 Pirogova str.}, \city{Novosibirsk}, \postcode{630090}, \country{Russia}}}


\abstract{{The paper presents the one of possible approach to model the epidemic propagation. The proposed model is based on the mean-field control inside separate groups of population, namely, suspectable (S), infected (I), removed (R) and cross-immune (C). In the paper the numerical algorithm to solve such a problem is presented, which ensures the conservation the total mass of population during timeline. Numerical experiments demonstrate the result of modelling   the propagation of COVID-19 virus during two 100 day periods in Novosibirsk (Russia).} }





\maketitle

\section{Introduction}\label{secIntrod}
The COVID-19 pandemic is not only accompanied by hundreds of thousands of deaths around the world, but also has a huge impact on the global economy, as well as on the daily life of everyone. As a rule, compartmental differential models of the SIR type are used to describe the dynamics of transmission of infectious diseases among the population. However, the SIR model stops working if it is necessary to take into account population heterogeneity or some stochastic phenomena that are significant in small populations and at the initial phase of the spread of the disease. Since the spread of COVID-19 has a significant spatial characteristic, various restrictive measures are used to slow down desease propogation, exerting an additional influence on the dynamics of population behavior  (by "spatial" we mean certain local features characteristic of a particular population, determined, for example, by location, politics, mass media, etc.). Thus, to study the spread of an infectious disease, it is also necessary to take into account the stochastic parameters of the system and its heterogeneity, which, in the classical sense, leads to systems with a large number of  differential equations, even if it is possible to cluster somehow the population. Computationally, the approach to describing population dynamics can be simplified using the ”Mean Field Game” theory, which allows to describe the behavioral dynamics of a system with a small number of partial differential equations, taking into account the heterogeneity of the population.

The paper organized as follows. In Section~\ref{secMatMod} SIRC model for disease spread is formulated based on system of ODEs and MFG approach is applied. The optimality conditions for mean-field system are obtained. In Section~\ref{secOptControl} the optimal control problem with external influence is obtained. In Section~\ref{secNumResults} the modification of computational scheme for MFG system is proposed and numerical experiments for COVID-19 disease spread in Russian Federation region are obtained and discussed.

\section{Mathematical model} \label{secMatMod}
There are a lot of mathematical models of infectious disease propagation based on mass balance law and described by systems of nonlinear ordinary differential equations~\cite{Brauer_2017}. The basic compartmental model to describe the transmission of communicable diseases is proposed by Kermack and McKendrick in 1927~\cite{KermackMcKendrick_1927} where the population is divided into three groups: susceptible, infected and removed. Various disease outbreaks, including the SARS epidemic of 2002-2003, the concern about a possible H5N1 influenza epidemic in 2005, the H1N1 influenza pandemic of 2009, the Ebola outbreak of 2014, and the COVID-19 pandemic in 2019 have re-ignited interest in epidemic models, beginning with the reformulation of the Kermack-McKendrick model. In the next section we applied one of variation of compartmental model where the population is divided into four groups taking into account cross-immune cases that take place while describing of the COVID-19 pandemic.

\subsection{The SIRC model} \label{secSIRCmod}
In this subsection, we will focus on the macroscopic description of the spread of infectious diseases without taking into account the individual affect of an agent on the whole population. For these purposes epidemiological models with different characteristics are generally used. Most of these models are based on the susceptible-infected-removed (SIR) separation. In work by Casagrandiet et al. \cite{Cas} a SIRC model is presented for describing the dynamic behavior of influenza A by introducing a new component, namely, cross-immunity (C) part of population for people who have recovered from infection with different strains of the same virus subtype at previous time intervals. In \cite{Zha, Mou, Sam}  this idea was developed and it was shown in \cite{Hai} that this model can be used to take into account the infected but asymptomatic part of the population.

\begin{figure}[!ht]
    \centering
    \includegraphics[width=0.7\textwidth]{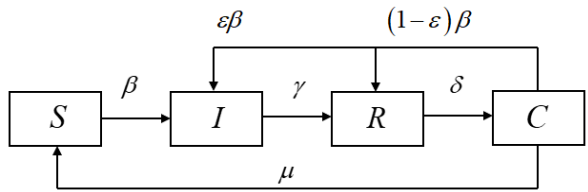}
    \caption{General flow diagram of the SIRC model}
    \label{figSIRCdiagramm}
\end{figure}

Consider the dynamic of the epidemic followed a simplified SIRC model. Assume that the epidemic occurs in a short time period in compare to the population dynamics (birth and death), that allows to neglect the last one. The flow diagram of the simplified SIRC model is presented in fig. \ref{figSIRCdiagramm}. From the SIRC point of view we divide the population at  any time moment $t$ into four compartments: the proportion of susceptible people $S(t)$, i.e., those who do not have specific immune defenses against the infection; the fraction $I(t)$ of those individuals who are infected; and two classes of those who are totally or partially immune  (i.e., recovered $R(t)$ and cross-immune $C(t)$, respectively).

Formally, the SIRC model in accordance with the diagram shown in fig.~\ref{figSIRCdiagramm} and mass balance law can be represented by the following system of four ordinary differential equations
\begin{equation}
\label{eq_1}
    \left\{ \begin{aligned}
  & \frac{dS\left( t \right)}{dt}=-\beta S\left( t \right)I\left( t \right)+\mu C\left( t \right), \\ 
 & \frac{dI\left( t \right)}{dt}=\beta S\left( t \right)I\left( t \right)+\varepsilon \beta C\left( t \right)I\left( t \right)-\gamma I\left( t \right), \\ 
 & \frac{dR\left( t \right)}{dt}=\left( 1-\varepsilon  \right)\beta C\left( t \right)I\left( t \right)+\gamma I\left( t \right)-\delta R\left( t \right), \\ 
 & \frac{dC\left( t \right)}{dt}=\delta R\left( t \right)-\beta C\left( t \right)I\left( t \right)-\mu C\left( t \right) \\ 
\end{aligned} \right.
\end{equation}
with initial conditions: $S(0)=S_0;\;I(0)=I_0;\;R(0)=R_0;\;C(0)=C_0$, where $S_0,\;I_0,\;R_0,\;C_0 \in \mathrm{R}^+$. Note, that $S(t)+I(t)+R(t)+C(t)$ represent the total population. Here the following parameters are used: $\mu$ is the rate at which the cross-immune population becomes susceptible again; $\beta$ is contact/infection transmission rate; $\varepsilon$ is the average reinfection probability of a cross-immune individual; $\gamma $ is recovery rate of the infected population; $\delta$ is the rate at which the recovered population becomes the cross-immune population and moves from total to partial immunity. Parameter $\mu$ is connected with average days when neutralizing antibodies disappeared, and $\delta$ is estimated via days when CD4+ and CD8+ T-cells decrease~\cite{Dan}. Note, that in the absence of cross-immunity, i.e. $(1-\varepsilon)=0$, the SIRC model becomes similar to the SIRS model since the two fractions $S$ and $C$ become immunologically indistinguishable \cite{Sam}. 

\subsection{Mean Field approach to epidemic control} \label{MFGmod}

In fact, greater consistency with some real-world phenomena can be achieved if stochastic changes in the system are taken into account. Epidemic models based on SIR approach assume that observed dynamics are determined by deterministic cases, making the assumption of population homogeneity, omitting the characteristics of individual agents. At the same time, taking into account the "rationality" of an individual who is part of population (that is, his ability to be in a state $X(t)$ and be able to change it) leads to management problems with a large number of participants in the system. Taking into account the stochasticity of the process, we assume that the dynamics of a single individual is subject to the $\hat{\text{I}}$to differential equation
\begin{equation}
\label{eq_2}
dX^N_i(t) = b (t,X^N_i(t),\theta^N(t),\alpha_i(t))dt + \sigma (t,X^N_i(t),\theta^N(t)) dW^N_i(t),   
\end{equation}
where $i \in {1,..,N};\; W^N_i $ is independent standart Wiener processes; $\alpha_i(t)$ is strategy of $i-$th agent and $\theta^N(t)$ is the empirical measure of distribution of agent (individuals, players) in the system at time $t$ \cite{Fis}. Put function $b,\; \sigma $ to be continuous in time and the same for all players. In \cite{Fis} it was shown that when the number of agents in system is extremely arising ($N \rightarrow \infty$), we can substitute the mass of single individuals to representative agent whose  state is determined by the following control equation 
\begin{equation}
\label{eq_3}
dX(t) = b (t,X(t),m(t),\alpha(t))dt + \sigma (t,X(t),m(t)) dW(t). 
\end{equation}
Here $X(t):[0,T]\rightarrow \Omega$; $m(t):\theta^N(t)\overset{N \rightarrow \infty}{\rightarrow} m(t)$ is agent's distribution over state space $\Omega$ in time $t$, and $\alpha(t)$ is the control process (or, in other words, strategy of representative agent), which ensures the Nash equilibrium of the system of interacting agents and minimize the function
\begin{equation}
\label{eq_4}
    J(\alpha) = \mathbf{E}\left[ \int_0^T f\left(s,X(s),m(t),\alpha(t)\right)ds + G \left(X(T),m(T)\right)\right], 
\end{equation}
where $f$ and $G$ are Lipschitz functions. Such an approach to get the control of a population with a large number of interacting agents is called the Mean Field Games (MFGs). For the more detailed description, we refer to the original work \cite {Fis, Ben, Las}, where the idea of MFGs originates.

For our goals consider the population of "atomized" agents characterized by their states (or positions) $x \in \Omega$ at each time moment $t \in [0,T]$. The term "atomized" means that each agent of an infinite set has no influence on the situation (because of its zero-measured support) and chooses the rational strategy $\alpha(t,x)$ taking into account its own position and the distribution $m(t,x): [0,T]\times \Omega \rightarrow \mathrm{R}$ of other agents. In \cite{Ben} it was shown that when the $\sigma$ in \eqref{eq_2},\eqref{eq_3}, which is characterized the stochastic nature of equilibrium of the process of agents interaction, is constant  the distribution of agents $m(t,x) $ obeys the Kolmogorov (Fokker-Planck) equation  
\begin{equation}
\label{eq_5}
\partial m/\partial t - \Delta m\sigma^2/2 + \nabla (m \alpha) = 0 \text{ in } [0,T]\times \Omega
\end{equation}
with initial condition 
\begin{equation}
\label{eq_6}    
m(0,x) = m_0(x) \text{ on } \Omega
\end{equation}
and Neumann boundary condition 
\begin{equation}
\label{eq_7}
\partial m/ \partial x = 0\; \forall t \text{ and } x \in \Gamma_\Omega.  
\end{equation}
Here boundary condition \eqref{eq_7} prevents the loss of density $m(t,x)$ with time. 

So, on the one hand, we have a differential SIRC model \eqref{eq_1} that describes the population at the macroscopic level, divide population on several clusters, but which is far from individual description. On the other hand, we can introduce "individuality" into the population and control its changing using the MFG model \eqref{eq_4}--\eqref{eq_7}. Following the ideas proposed in \cite{Lee} we combine the  spatial SIRC model and MFG and introduce a management problem to control the virus spreading within several groups of population. 

Instead of $S,I,R,C$ denoted the masses of the groups of population we introduce the density of distribution of agents within these groups $m_i(t,x):[0,T]\times[0,1]\rightarrow \mathrm{\textbf{R}}$, where $i \in\{S,I,R,C\}$. Put that the state variable $x$ varies within $\Omega = [0,1]$ and denotes the population's propensity to comply with quarantine measures (in particular, physical distancing and self-isolation). Here $0$ means that the agent is loyal to the introduced restriction measures and inclined to comply with them, and $1$ means exactly the opposite. Introduce also the functions $\alpha_i(t,x):[0,T]\times[0,1]\rightarrow \mathrm{\textbf{R}}$, $i \in\{S,I,R,C\}$ which denote the strategy of compliance of the self-isolation of representative agent in each group of population. Due to \eqref{eq_1}, \eqref{eq_5} -- \eqref{eq_7} the functions $m_i(t,x)$ are the subject of the following PDEs system:
\begin{equation}
\label{eq_8}
 \hspace{-4mm}\left\{ \begin{aligned}
   & {{\partial }_{t}}{{m}_{S}}+\nabla \left( {{m}_{S}}{{\alpha }_{S}} \right)+\beta {{m}_{S}}{{m}_{I}}-\mu {{m}_{C}}-{\sigma _{S}^{2}\Delta {{m}_{S}}}/{2}\;=0, \\ 
  & {{\partial }_{t}}{{m}_{I}}+\nabla \left( {{m}_{I}}{{\alpha }_{I}} \right)-\beta {{m}_{S}}{{m}_{I}}-\varepsilon \beta {{m}_{C}}{{m}_{I}}+\gamma {{m}_{I}}-{\sigma _{I}^{2}\Delta {{m}_{I}}}/{2}\;=0, \\ 
 & {{\partial }_{t}}{{m}_{R}}+\nabla \left( {{m}_{R}}{{\alpha }_{R}} \right)-(1-\varepsilon )\beta {{m}_{C}}{{m}_{I}}-\gamma {{m}_{I}}+\delta {{m}_{R}}-{\sigma _{R}^{2}\Delta {{m}_{R}}}/{2}\;=0, \\ 
  & {{\partial }_{t}}{{m}_{C}}+\nabla \left( {{m}_{C}}{{\alpha }_{C}} \right)-\delta {{m}_{R}}+\beta {{m}_{C}}{{m}_{I}}+\mu {{m}_{C}}-{\sigma _{C}^{2}\Delta {{m}_{C}}}/{2}\;=0. \\ 
\end{aligned} \right.
\end{equation}
Here $\sigma_i$, $i \in \{S,I,R,C\}$ are non-negative parameters which as before characterize the stochastic processes within the population. Put that the initial value of $m_i$ are given
\begin{equation}
\label{eq_9}
 m_i (0,x)=m_{i0}(x)\;\; \forall x\in[0,1].
\end{equation}
In addition, similarly with \eqref{eq_7} put zero Neumann boundary conditions for $m_i$, that is no mass can flow in or out of $\Omega$.
\begin{equation}
\label{eq_10}
\partial m_i/ \partial x = 0\; \;\forall t \in [0,T] \;\text{ and } \;x = 0,1.  
\end{equation}

Note that  
\begin{equation}
\nonumber
\displaystyle\int_{0}^{1} \left(m_S(t,x) + m_I(t,x) + m_R(t,x) + m_C(t,x) \right) \text{d}x = \int_{0}^{1} m(t,x) \text{d}x, \; \forall t\in [0,T],
\end{equation}
where $m(t,x)$ is total mass of population. Also note that summing up the equations in system \eqref{eq_8}, we obtain equation \eqref{eq_5} with respect to the total mass of population $m(t,x)$. It leads to the following requirement
\begin{equation}
\nonumber
\frac{\partial}{\partial t} \displaystyle\int_{0}^{1} \left(m_S(t,x) + m_I(t,x) + m_R(t,x) + m_C(t,x) \right) dx = 0,
\end{equation}
i.e., the total mass of the four groups will be conserved $\forall t$.

Now, based on the principles of MFG, assume that agents are rational and tend to choose the strategy to maximize their own benefit in accordance with cost functions of the form of \eqref{eq_4}. Produce the following management problem: \textit{minimize the value functional in respect to $(m_i,\alpha_i)\; \forall i \in \{S,I,R,C\}$}
\begin{multline}
\label{eq_11}
J (m_{SIRC},\alpha_{SIRC})= \int_{0}^{T} {\int_{0}^{1}\sum_{i\in\{S,I,R,C\}}  
 ( F_i\left( \alpha_{i} ,t,x \right)m_i + } \\
 +g_i(t,x,m_{i}) )dxdt + \int_{0}^{1} m_I^2(T,x)/2\, dx,
\end{multline}
where index ${SIRC}$ denotes the set of all possible value from $\{S,I,R,C\}$. Here $F_i$ function denotes the running cost of implementation of strategy $\alpha_{i} $; $g_i$ is payment for current position or state and $\int_{0}^{1} m_I(T,x) dx$ is the terminal cost for epidemic propagation. 

As the function of cost of the strategy implementation $F_i(\bar{\alpha}_{i},t,x)$ (the agent's payment for state changing), consider a piecewise continuous function
\begin{equation}
\label{eq_12}
F_i(\bar{\alpha}_{sirc} ,t,x)=\left\{\begin{array}{l}
\bar{\alpha }_i^{2}(100+\bar{\alpha_i}^2) / 2 \;\;\; {\rm if }\, \,\,\,\bar{\alpha_i }\le 0,\\
\bar{\alpha }_i^{2}(20+\bar{\alpha_i}^2) / 2 \,\,\,\,\,\; {\rm otherwise} 
\end{array}\right.
\end{equation}
for all admissible values of $\bar{\alpha_i} \in \mathrm{\textbf{R}}$ and $\forall i \in \{S,I,R,C\}$. The form of \eqref{eq_12} means that the transition to self-isolation for agent is more expensive than maintaining close contact with other people. Note, that for \eqref{eq_12} the following properties are satisfied $\forall \left( t,x \right)\in [0,T]\times [0,1]$:
\begin{equation}
 \label{eq_13}
\begin{aligned}
& \hspace{-25mm} 1.\;		\frac{\partial F_i}{\partial \bar{\alpha_i }}\left( \bar{\alpha_i },t,x \right) \text{ is continuous for all admissible } \bar{\alpha_i }; \\
& \hspace{-25mm}2. \; 	F_i\left( 0,t,x \right)=\frac{\partial F_i}{\partial \bar{\alpha }}\left( 0,t,x \right)=0; \;\\
& \hspace{-25mm} 3.\;  \frac{\partial F_i}{\partial \bar{\alpha_i }}\left( \bar{\alpha_i },t,x \right) \text{ is strictly monotonous for } \bar{\alpha_i }\in \left( -\infty ,+\infty  \right).
\end{aligned}  
\end{equation}
These properties ensure the unique solvability of the equation 
\[{\partial F_i}/{\partial \bar{\alpha_i }}\;\left( \bar{\alpha_i },t,x \right)=z \text{ for any } z\in \mathrm{\textbf{R}} \text{ and } \left( t,x \right)\in [0,T]\times [0,1].\]

For $g_i(t,x,m_{i})$ function we put 
\begin{equation}
\label{eq_14}
g_i(t,x,m_{i})=\left\{\begin{array}{l}
\displaystyle\frac{c_{1i} (1-x)}{1 +c_{2i} m_i}  m_i \;\;\; \text{ if } i \in \{S,R,C\}, \vspace{2mm}\\
c_{1I}m_Ix, \,\,\,\,\,\;\text{ if } i =I.
\end{array}\right.
\end{equation}
Here $c_{1i}, \; c_{2i} \in [0,1]$ are positive constants. This concept means that if a
person is not infected and the number of non-isolating people in population is arising, person's profit decreases
to comply with the restrictions. This effect is commonly called "economies of scale" \cite{McC}. On the contrary, we believe that if a person is infected, then he is inclined to comply with self-isolation.

Thus, the MFG approach generates the following optimization problem: find the optimal strategy $\alpha_i$ by minimizing the function \eqref{eq_11} with the restriction  in the form of PDEs system \eqref{eq_8} with initial \eqref{eq_9} and boundary \eqref{eq_10} conditions. 

\subsection{Optimality conditions} \label{secOptCond}

To ensure that the strategy chosen by a person is optimal use the Lagrange multiplier method \cite{Ben}. Multiply the first equation in \eqref{eq_8} by an arbitrary smooth function $ v_S(t, x)\in C ^ {\infty} \left(\left[0, T \right] \times [0,1] \right) $ and integrate the resulting expression by parts with respect to $ t $ and $ x $:
\begin{equation}
 \begin{aligned}
\label{eq_15}
&{L_S}:=-\int_{0}^{T}{\int_{0}^{{1}}{\left( {\partial v_S}/{\partial t}\;+{{\sigma_S }^{2}}\Delta v_S/2+\alpha_S \cdot  \partial v_S/\partial x \right)m_S\,\text{d}x \,}\text{d}t} + \\
&\hspace{20mm}+\int_{0}^{T}{\int_{0}^{{1}}}({\beta {{m}_{S}}{{m}_{I}}v_S-\mu {{m}_{C}}v_S})\text{d}x\text{d}t + \\
&\hspace{30mm}+\int_{0 }^{{1}}{\left( v_S(T,x)m_S(T,x)-v_S(0,x){{m}_{S0}}(x) \right)}\text{ d}x =0. 
\end{aligned}   
\end{equation}
Do the same with other equations in \eqref{eq_8} for smooth function $ v_i(t, x)\in C ^ {\infty} \left(\left[0, T \right] \times [0,1] \right) $ where ${i \in \{I,R,C\}}$:
\begin{equation}
 \begin{aligned}
\label{eq_16}
&{L_I}:=-\int_{0}^{T}{\int_{0}^{{1}}{\left( {\partial v_I}/{\partial t}\;+{{\sigma_I }^{2}}\Delta v_I/2+\alpha_I \cdot  \partial v_I/\partial x \right)m_I\,\text{d}x \,}\text{d}t} + \\
&\hspace{20mm}+\int_{0}^{T}{\int_{0}^{{1}}}(-\beta {{m}_{S}}{{m}_{I}}v_I-\varepsilon \beta {{m}_{C}}{{m}_{I}}v_I+\gamma {{m}_{I}}v_I)\text{d}x\text{d}t + \\
&\hspace{30mm}+\int_{0 }^{{1}}{\left( v_I(T,x)m_I(T,x)-v_I(0,x){{m}_{I0}}(x) \right)}\text{ d}x =0. 
\end{aligned}   
\end{equation}
\begin{equation}
 \begin{aligned}
\label{eq_17}
&{L_R}:=-\int_{0}^{T}{\int_{0}^{{1}}{\left( {\partial v_R}/{\partial t}\;+{{\sigma_R }^{2}}\Delta v_R/2+\alpha_R \cdot  \partial v_R/\partial x \right)m_R\,\text{d}x \,}\text{d}t} + \\
&\hspace{15mm}+\int_{0}^{T}{\int_{0}^{{1}}}(-(1-\varepsilon )\beta {{m}_{C}}{{m}_{I}}v_R-\gamma {{m}_{I}}v_R+ \delta {{m}_{R}}v_R)\text{d}x \text{d}t  +\\
&\hspace{25mm}+\int_{0 }^{{1}}{\left( v_R(T,x)m_R(T,x)-v_R(0,x){{m}_{R0}}(x) \right)}\text{ d}x =0. 
\end{aligned}   
\end{equation}
\begin{equation}
 \begin{aligned}
\label{eq_18}
&{L_C}:=-\int_{0}^{T}{\int_{0}^{{1}}{\left( {\partial v_C}/{\partial t}\;+{{\sigma_C }^{2}}\Delta v_C/2+\alpha_C \cdot  \partial v_C/\partial x \right)m_C\,\text{d}x \,}\text{d}t} + \\
&\hspace{20mm}+\int_{0}^{T}{\int_{0}^{{1}}}(-\delta {{m}_{R}}v_C+\beta {{m}_{C}}{{m}_{I}}v_C+\mu {{m}_{C}}v_C)\text{d}x \text{d}t  +\\
&\hspace{25mm}+\int_{0 }^{{1}}{\left( v_C(T,x)m_C(T,x)-v_C(0,x){{m}_{C0}}(x) \right)}\text{ d}x =0. 
\end{aligned}   
\end{equation}
Note that expressions \eqref{eq_15}--\eqref{eq_18} are valid when the following boundary conditions are satisfied
\begin{equation}
\label{eq_19}
\partial v_i/ \partial x = 0\; \;\forall t \in [0,T] \;\text{ and } \;x = 0,1\; {\forall i \in \{S,I,R,C\}}  
\end{equation}
and 
\begin{equation}
\label{eq_20}
\alpha_i(t,0) = \alpha_i(t,1) = 0\; \;\forall t \in [0,T] \; {\forall i \in \{S,I,R,C\}}. 
\end{equation}
Now write down the Lagrange function corresponding to the optimization problem under consideration
\begin{equation}
\begin{aligned}
\label{eq_21}
\Im (m_{SIRC},\alpha_{SIRC},v_{SIRC}):=J (m_{SIRC},\alpha_{SIRC}) - L_S - L_I - L_R - L_C.
\end{aligned}    
\end{equation}
As the result, the minimization of (\ref{eq_11}) together with the  (\ref{eq_8})--(\ref{eq_10}) can be represented \cite{Ben} as the problem of finding a saddle point
\begin{equation}
\label{eq_22}
\underset{(m_i,\alpha_i )}{\mathop{\inf }}\,\,\,\underset{v_i}{\mathop{\sup}}\,\,\,\Im (m_{SIRC},\alpha_{SIRC},v_{SIRC}) \;{\forall i \in \{S,I,R,C\}}.  
\end{equation} 
Variation (\ref{eq_15})--(\ref{eq_18}) with respect to components $ m_i\; \forall (t,x)\in [0,T]\times[0,1]$  to find a stationary point gives the system of analogues of Hamilton-Jacobi-Bellman equation
\begin{equation}
\label{eq_23}
 \hspace{-5mm}\left\{ \begin{aligned}
   & {\partial v_S}/{\partial t}\;+{{\sigma_S }^{2}}\Delta v_S/2+\alpha_S \cdot \partial v_S / \partial x + \beta m_I (v_I-v_S) =-F_S\,-{\partial g_S}/{\partial m_S}, \\ 
  & \begin{aligned}
  {\partial v_I}/{\partial t}\;&+{{\sigma_I }^{2}}\Delta v_I/2+\alpha_I \cdot \partial v_I / \partial x  + \beta m_S (v_I-v_S) + \beta m_C (v_R-v_C)  \\
  &+ (\varepsilon \beta m_C - \gamma)(v_I-v_R)  =-F_I\,-{\partial g_I}/{\partial m_I} -\delta(T-t)m_I,
  \end{aligned} \\ 
 & {\partial v_R}/{\partial t}\;+{{\sigma_R }^{2}}\Delta v_R/2+\alpha_R \cdot \partial v_R / \partial x + \delta (v_C-v_R) =-F_R\,-{\partial g_R}/{\partial m_R} , \\ 
  &  \begin{aligned}
  {\partial v_C}/{\partial t}\;&+{{\sigma_C }^{2}}\Delta v_C/2+\alpha_C \cdot \partial v_C / \partial x + \mu (v_S-v_C) + \\
  &+ \varepsilon \beta m_I (v_I-v_R) + \beta m_I (v_R - v_C)  =-F_C\,-{\partial g_C}/{\partial m_C}, 
  \end{aligned} \\ 
\end{aligned} \right.
\end{equation}
where $\delta(T-t)$ is delta function and zero initial condition 
\begin{equation}
\label{eq_24}  
v_i(T,\,x)=0\,\,\,\,\,\forall \text{ }x\in \left[ 0,1 \right].
\end{equation}
Variation (\ref{eq_15})--(\ref{eq_18}) with respect to components $ \alpha_i\; \forall (t,x)\in [0,T]\times[0,1]$   gives the following optimality conditions for $\bar{\alpha_i} \in \mathrm{R}$  in addition to system \eqref{eq_23}, \eqref{eq_24}, \eqref{eq_19}
\begin{equation}
\label{eq_25}  
\frac{\partial F_i}{\partial \bar{\alpha_i }}\left( \alpha_i ,t,x \right)+\frac{\partial v_i}{\partial x}\left( t,x \right)=0
\end{equation}
$\forall\, i\in \{S,I,R,C\}\,\forall \left( t,x \right)\in [0,T]\times \left[ 0,1 \right]$.
The properties \eqref{eq_13} ensure that \eqref{eq_25} has a unique solution that satisfies \eqref{eq_20}.

Thus, two coupled systems of PDEs with initial and boundary conditions \eqref{eq_8} -- \eqref{eq_10} and \eqref{eq_19}, \eqref{eq_23}, \eqref{eq_24} together with \eqref{eq_25} give the necessary conditions for the minimization of \eqref{eq_11}.

\section{Optimal control problem with external influence}\label{secOptControl}

In our case, we consider agents "rational" but "selfish", i.e. inclined to choose a strategy that maximizes their own benefits, based on the strategies of other agents. This can lead to a faster spread of the infection. In real case, the situation is controlled by some external influence (the introduction of strict quarantine and other restrictions or the allocation of social benefits, etc.). Such quarantine norms become decisive in the behavior of agents and, as a rule, are strict requirements, therefore, they cannot be taken into account in the cost functional, where their non-fulfillment becomes possible. Thus, assume that the total control of behavior of system of agents can be determined by two processes:
\begin{equation}
\label{eq_26}
 \alpha_i = \hat{\alpha_i }\left( t,x \right)+\rho_i \left( \hat{\alpha_i }\left( t,x \right),t,x \right).  \\
\end{equation}
Here $\hat{\alpha_i }$, as before, is the person's strategy and $\rho_i (\hat{\alpha_i },t,x)$ is the external (corrective) control which adjusts the player’s strategy depending on the fulfillment of a certain condition. For example, we can put the following expression: 
\begin{equation}
\label{eq_27}
\rho_i \left( \hat{\alpha_i },t,x \right)=\left\{ \begin{matrix}
    0 \,\,\,\,\,\,\;\;\;\;\;\;\;\,\,\,\text{if}\,\,\hat\alpha_i\le 0,  \\
   - c_{3i}\hat{\alpha_i }\,\,\,\,\text{otherwise},  \\
\end{matrix} \right.
\end{equation}
where $c_{3i} \in [1,2] $ are positive constants and \eqref{eq_27} does not violate the boundary condition \eqref{eq_20} for $\alpha_i$. From physical point of view \eqref{eq_27} means that if agents are inclined to observe the quarantine measures, the government doesn't impose additional restrictions, otherwise takes some restriction actions. As before, we define the dynamics of the agent's mass  as the solution of the system of equations \eqref{eq_8} taking into account representation \eqref{eq_26}. Instead of \eqref{eq_11} we will minimize the same function $J( m_{SIRC},\hat\alpha_{SIRC})$ but for set of ${\hat{\alpha_i}}\; \forall i \in \{S,I,R,C\}$ 
\begin{multline}
\label{eq_28}
J (m_{SIRC},\hat\alpha_{SIRC})= \int_{0}^{T} {\int_{0}^{1}\sum_{i\in\{S,I,R,C\}}  
 ( F_i\left( \hat\alpha_{i} ,t,x \right)m_i + } \\
 +g_i(t,x,m_{i}) )dxdt + \int_{0}^{1} m_I^2(T,x)/2\, dx.
\end{multline}
Thus, we do not take into account the cost of implementing external control, because it is not important for a single individual. This assumption gives us the optimality conditions in the form of the system \eqref{eq_23} with \eqref{eq_26} and instead of \eqref{eq_25} we get 
\begin{equation}
\label{eq_29}
\frac{\partial F_i}{\partial \bar{\hat{\alpha_i }}}\;+\left( 1+\frac{{\partial\rho_i \left( {\hat{\alpha_i },t,x} \right)}}{{\partial\bar{\hat{\alpha_i }}}}\; \right)\frac{\partial v_i}{\partial x}=0\;\,\,\,\text{in}\,\,\,[0,T]\,\times \,[0,1].  
\end{equation}
Thus,  the  following  optimization  problem with external control can be formulated: in accordance with \eqref{eq_26} find the minima of the function \eqref{eq_28} with the restriction \eqref{eq_8} -- \eqref{eq_10},  \eqref{eq_19}, \eqref{eq_23}, \eqref{eq_24} and \eqref{eq_29}. 

\section{Numerical experiments}\label{secNumResults}
We provide the numerical experiments to show a difference between standard differential SIRC model \eqref{eq_1} and its MFG extension. 

\subsection{Numerical scheme for FPK equations}

To search for the numerical implementation of differential problem mentioned above, introduce discrete uniform
grids in time
$$
t_{k} =k\tau ,\, \, \, \, \, k=0,...,M,\, \, \, \tau =T\big/ M;
$$
and in space
$$
x_{i+1/2} =(i+1/2)h,\, \, \, \,\, i=-1,...,N,\, \, \, h={1/N} .
$$

We will serch for a solution as a piecewise linear functions $m^{i,h} (t,x)$ at each time level $t_{k} $, which are continuous on [0,1] and linear on each segment $\omega _{j} =\left[x_{j-1/2} ,x_{j+1/2} \right]\,\,\forall \, j=1,...,N-1$. In addition, $m^{i,h} (t_{k} ,x)$ is assumed to be constant on the intervals $\omega _{0} =[0,x_{1/2}]$ and $\omega _{N} =\left[x_{N-1/2} ,1\right]$ to satisfy boundary condition \eqref{eq_10}. Thus, at each time level $t_{k} $, the function $m^{i,h} (t_{k} ,x)$ is completely determined by its discrete values $m_{k,j+{1/2} }^{i,h} := m^{i,h} \left(t_{k} ,x_{j+{1/2} } \right):$ 
$$ m^{i,h} (t_{k} ,x)=m_{k,j-{1/2} }^{i,h} {(x_{j+{1/2} } -x)\big/ h} +m_{k,j {1/2} }^{h} {(x-x_{j-{1/2} } ) \big/ h} \, \, \, \, \, \, \forall \, x\in \omega _{i} . $$
To satisfy the boundary condition \eqref{eq_10}, we put
\begin{equation}
\label{eq_30}
m_{k,-{1/2} }^{i,h} =m_{k,{1/2} }^{i,h} \,\,\, {\rm and} \,\,\,\,
m_{k,N+{1/2} }^{i,h} =m_{k,N-{1/2} }^{i,h} .
\end{equation}

To construct a computational scheme for equations \eqref{eq_8}, modify the idea firstly proposed by  prof. Shaydurov V.V. for numerical solution of FPK equation \eqref{eq_5} and applied to various MFG statements in \cite{Sha1,Sha2,Sha3,Sha4}. The idea was to split the approximation of the FPK \eqref{eq_5} equation into two parts: advection and diffusion one. Note, that advection part of equations in \eqref{eq_8} can be written in the following form
\begin{equation}
 \label{eq_31}
 \frac{\partial {{m}_{i}} }{\partial t}+\frac { \partial \left( {{m}_{i}}{{\alpha }_{i}} \right)}{\partial x} = f^i_1(t,x),    
\end{equation}
$\forall i\in\{S,I,R,C\} $ where $f^i_1(t,x)$ are
\begin{equation}
\label{eq_32}
 \hspace{-5mm}\left\{ \begin{aligned}
   & f^S_1 = -\beta {{m}_{S}}{{m}_{I}}+\mu {{m}_{C}}+{\sigma _{S}^{2}\Delta {{m}_{S}}}/{2}\;, \\ 
  & f^I_1 = \beta {{m}_{S}}{{m}_{I}}+\varepsilon \beta {{m}_{C}}{{m}_{I}}-\gamma {{m}_{I}}+{\sigma _{I}^{2}\Delta {{m}_{I}}}/{2}\;, \\ 
 & f^R_1 = (1-\varepsilon )\beta {{m}_{C}}{{m}_{I}}+\gamma {{m}_{I}}-\delta {{m}_{R}}+{\sigma _{R}^{2}\Delta {{m}_{R}}}/{2}\;, \\ 
  & f^C_1 = \delta {{m}_{R}}-\beta {{m}_{C}}{{m}_{I}}-\mu {{m}_{C}}+{\sigma _{C}^{2}\Delta {{m}_{C}}}/{2}\;=0. \\ 
\end{aligned} \right.
\end{equation}
For diffusion part of equations in \eqref{eq_8} use the same property 
\begin{equation}
\label{eq_33}
 \hspace{-5mm}\left\{ \begin{aligned}
   & \beta {{m}_{S}}{{m}_{I}}-\mu {{m}_{C}}-{\sigma _{S}^{2}\Delta {{m}_{S}}}/{2}\;= f^S_2 (t,x), \\ 
  & -\beta {{m}_{S}}{{m}_{I}}-\varepsilon \beta {{m}_{C}}{{m}_{I}}+\gamma {{m}_{I}}-{\sigma _{I}^{2}\Delta {{m}_{I}}}/{2}\;=f^I_2 (t,x), \\ 
 & -(1-\varepsilon )\beta {{m}_{C}}{{m}_{I}}-\gamma {{m}_{I}}+\delta {{m}_{R}}-{\sigma _{R}^{2}\Delta {{m}_{R}}}/{2}\;=f^R_2 (t,x), \\ 
  & -\delta {{m}_{R}}+\beta {{m}_{C}}{{m}_{I}}+\mu {{m}_{C}}-{\sigma _{C}^{2}\Delta {{m}_{C}}}/{2}\;=f^C_2 (t,x), \\ 
\end{aligned} \right.
\end{equation}
where $f^2_i = - \partial_t m_i  - \nabla (m_i \alpha_i).$  Following the idea proposed in \cite{Sha1,Sha2,Sha3,Sha4} get the approximation for \eqref{eq_31} in the point $\left( {{t}_{k,}}{{x}_{i+1/2}} \right)$ as
\begin{equation}
\label{eq_34}
\begin{aligned}
&\frac{1}{8\tau }m_{k,j-{1}/{2}\;}^{i,h}+\frac{3}{4\tau }m_{k,j+{1}/{2}\;}^{i,h}+\frac{1}{8\tau }m_{k,j+{3}/{2}\;}^{i,h} - \gamma _{k,j+{1}/{2}\;}^{i,1}m_{k-1,j-{1}/{2}\;}^{i,h}- \\
&-\gamma _{k,j+{1}/{2}\;}^{i,2}m_{k-1,j+{1}/{2}\;}^{i,h}-\gamma _{k,j+{1}/{2}\;}^{3}m_{k-1,j+{3}/{2}\;}^{i,h} = f^i_1\left( {{t}_{k,}}{{x}_{j+1/2}} \right)
\end{aligned}
\end{equation}
$\forall i\in\{S,I,R,C\} $ and $\forall \,j=0,...,N-1\,\,\text{ }\forall \,k=1,...,M$ with coefficients
\begin{equation}
\label{eq_35}
\begin{aligned}
&\gamma _{k,j+{1}/{2}\;}^{i,1}=\frac{1}{8\tau }\left( 1+\frac{4\tau }{h}{{\alpha}^i_{k,j}} \right),\,\,\,\\
&\gamma _{k,j+{1}/{2}\;}^{i,2}=\frac{1}{8\tau }\left( 3+\frac{4\tau }{h}{{\alpha }^i_{k,j}} \right)+\frac{1}{8\tau }\left( 3-\frac{4\tau }{h}{{\alpha ^i}_{k,j+1}} \right),\,\,\,\\
&\gamma _{k,j+{1}/{2}\;}^{i,3}=\frac{1}{8\tau }\left( 1-\frac{4\tau }{h}{{\alpha }^i_{k,j+1}} \right).\\
\end{aligned}
\end{equation}
Introduce the following notation for the approximation of the differential expression $\sigma^2_i \Delta m_i/2$ $\forall i \in \{S,I,R,C\}$:
\begin{equation}
\label{eq_37}
\Delta^{i,h}_{k,j+1/2}={\sigma _{i}^{2}}\frac{ m^{i,h}_{k,j-1/2} - 2 m^{i,h}_{k,j+1/2} + m^{i,h}_{k,j+3/2}   }{2h^2}.
\end{equation}

For approximation of \eqref{eq_33} use the following finite-difference scheme 
\begin{equation}
\label{eq_36}
 \hspace{-5mm}\left\{ \begin{aligned}
 & \beta m^{S,h}_{k-1,j+1/2}  m^{I,h}_{k-1,j+1/2}- \mu m^{C,h}_{k-1,j+1/2} - \Delta^{S,h}_{k,j+1/2} \approx f^S_2(t_k,x_{j+1/2}),\\
  &\begin{aligned}
  -\beta {{m}^{S,h}_{k-1,j+1/2}}{m}^{I,h}_{k-1,j+1/2}-\varepsilon\beta {{m}^{C,h}_{k-1, j+1/2}}&{{m}^{I,h}_{k-1, j+1/2}} + \gamma{{m}^{I,h}_{k-1, j+1/2}} -  \\
  &- \Delta^{I,h}_{k,j+1/2} \approx f^I_2(t_k,x_{j+1/2}),
  \end{aligned} \\
  &\begin{aligned}
    -(1-\varepsilon )\beta  {{m}^{C,h}_{k-1,j+1/2}}{{m}^{I,h}_{k-1,j+1/2}}-\gamma &{{m}^{I,h}_{k-1,j+1/2}}+\delta {{m}^{R,h}_{k-1, j+1/2}}- \\
    &-\Delta^{R,h}_{k,j+1/2} \approx f^R_2(t_k,x_{j+1/2}),
  \end{aligned}\\
  &\begin{aligned}
   -\delta {{m}^{R,h}_{k-1,j+1/2}}  + \beta {{m}^{C,h}_{k-1,j+1/2}}{{m}^{I,h}_{k-1,j+1/2}} & + \mu {{m}^{C,h}_{k-1,j+1/2}} - \\
   &-\Delta^{C,h}_{k,j+1/2} \approx f^C_2(t_k,x_{j+1/2}).
  \end{aligned} \\
\end{aligned} \right.
\end{equation}
Note that 
\[f_1^i(t_k,x_{j+1/2})+f_2^i(t_k,x_{j+1/2}) = 0 \;\;\]
$\forall k=1,..,M,\,j=1,...,N-1 \;\forall i\in \{S,I,R,C\}.$ It leads to the following finite-difference approximation of equations \eqref{eq_8}
\begin{equation}
\label{eq_39}
 D^{i,h}_{k,j+1/2} =  f^{i,h}_{k-1,j+1/2} + \Gamma^{i,h}_{k-1,j+1/2},   
\end{equation}
where 
\begin{equation}
\label{eq_40}
\begin{aligned}
&
\begin{aligned}
 D^{i,h}_{k,j+1/2} = \left(\frac{1}{8\tau }-\frac{\sigma^2_i}{2h^2}\right)m_{k,j-{1}/{2}\;}^{i,h}&+\left(\frac{3}{4\tau }+\frac{\sigma^2_i}{h^2}\right)m_{k,j+{1}/{2}\;}^{i,h} +\\
 &+\left(\frac{1}{8\tau }-\frac{\sigma^2_i}{2h^2}\right)m_{k,j+{3}/{2}\;}^{i,h},  
\end{aligned}
\\
&\begin{aligned}
\Gamma^i_{k-1, j+1/2} = \gamma _{k,j+{1}/{2}\;}^{i,1}m_{k-1,j-{1}/{2}\;}^{i,h} + 
 \gamma _{k,j+{1}/{2}\;}^{i,2}&m_{k-1,j+{1}/{2}\;}^{i,h}  + \\
& +\gamma _{k,j+{1}/{2}\;}^{i,3}m_{k-1,j+{3}/{2}\;}^{i,h}   
\end{aligned}
\end{aligned}
\end{equation}
and 
\begin{equation}
\label{eq_41}
 \left\{ \begin{aligned}
 & \begin{aligned}
  f^{S,h}_{k-1,j+1/2} = &-\beta m^{S,h}_{k-1,j+1/2}  m^{I,h}_{k-1,j+1/2} + \mu m^{C,h}_{k-1,j+1/2},
 \end{aligned}\\
  & \begin{aligned}
  f^{I,h}_{k-1,j+1/2} =\, \beta &{{m}^{S,h}_{k-1,j+1/2}}{m}^{I,h}_{k-1,j+1/2}+\\ &+\varepsilon\beta {{m}^{C,h}_{k-1, j+1/2}}{{m}^{I,h}_{k-1, j+1/2}} -\gamma{{m}^{I,h}_{k-1, j+1/2}},
 \end{aligned}\\
 & \begin{aligned}
  f^{R,h}_{k-1,j+1/2} =\,   (1-\varepsilon )\beta  {{m}^{C,h}_{k-1,j+1/2}}&{{m}^{I,h}_{k-1,j+1/2}}+\\
  &+\gamma {{m}^{I,h}_{k-1,j+1/2}}-\delta {{m}^{R,h}_{k-1, j+1/2}}, \\
 \end{aligned}\\
& \begin{aligned}
f^{C,h}_{k-1,j+1/2} =\, &\delta {{m}^{R,h}_{k-1,j+1/2}}  - \beta {{m}^{C,h}_{k-1,j+1/2}}{{m}^{I,h}_{k-1,j+1/2}}  - \mu {{m}^{C,h}_{k-1,j+1/2}},\\
\end{aligned}
\end{aligned} \right.
\end{equation}
and initial conditions corresponding to (\ref{eq_9})
\begin{equation}
\label{eq_FPK_init}
m_{0,j+{1}/{2}\;}^{i,h}={{m}^i_{0}}({{x}_{j+{1}/{2}\;}})\,\,\forall \,\,j=0,...,N-1.
\end{equation}

Thus, instead of differential equations \eqref{eq_8} $\forall i\in \{\S,I,R,C\}$ we get the following systems of algebraic equation 
\begin{equation}
\label{eq_42}
{\it {\mathfrak A^i}}m_{\cdot ,\cdot }^{i,h} =\left[\begin{array}{cccc}
{{\it {\mathbb A}}} & {} & {} & {} \\ {-{\it {\mathbb B}}^i_{2} } & {{\it {\mathbb A}}} & {} & {} \\
{} & {\ddots } & {\ddots } & {} \\ {} & {} & {-{\it {\mathbb B}}^i_{M} } & {{\it {\mathbb A}}} \end{array}\right]
\left[\begin{array}{c} {m_{1,\cdot }^{h} } \\ {m_{2,\cdot }^{h} } \\ {\vdots } \\
{m_{M,\cdot }^{h} } \end{array}\right]=\left[\begin{array}{c} {{\it {\mathbb B}}^i_{1} m_{0,\cdot }^{h}  + f^i_{0,\cdot}} \\
{f^i_{1,\cdot}} \\ {\vdots } \\ { f^i_{M-1,\cdot}} \end{array}\right]={\it {\mathfrak F^i}}m_{0,\cdot }^{h} .
\end{equation}
Here and later we use the points in subscript
to indicate acceptable values in corresponding position of mesh function.
${\it {\mathbb A}}$
and ${\it {\mathbb B}}^i_{k} $ are the matrices of the form
$$
{\it {\mathbb B}}^i_{k} =\left[\begin{array}{cccc} {\gamma _{k,-{1/2} }^{i,1} +
\gamma _{k,{1/2} }^{i,2} } & {\gamma _{k,{3/2} }^{i,3} } & {} & {} \\
{\gamma _{k,{1/2} }^{i,1} } & {\gamma _{k,{3/2} }^{i,2} } & {\gamma _{k,{5/2} }^{i,3} } & {} \\
{} & {\ddots } & {\ddots } & {} \\ {} & {} & {\gamma _{k,{N-3/2} }^{i,1} } & {\gamma _{k,N-{1/2} }^{i,2} +
\gamma _{k,N+{1/2} }^{i,3} }
\end{array}\right],
$$
$$
{\it {\mathbb A}}=\left[\begin{array}{cccc} {\displaystyle\frac{7}{8\tau } +
\displaystyle\frac{\sigma ^{2} }{2h^{2} } } & {\displaystyle\frac{1}{8\tau } -\frac{\sigma ^{2} }{2h^{2} } } & {} & {} \\ \\
{\displaystyle\frac{1}{8\tau } -\frac{\sigma ^{2} }{2h^{2} } } & {\displaystyle\frac{3}{4\tau } +\frac{\sigma ^{2} }{h^{2} } } &
{\displaystyle\frac{1}{8\tau } -\frac{\sigma ^{2} }{2h^{2} } } & {} \\ {} & {\ddots } & {\ddots } & {} \\
{} & {} & {\displaystyle\frac{1}{8\tau } -\frac{\sigma ^{2} }{2h^{2} } } & {\displaystyle\frac{7}{8\tau } +\frac{\sigma ^{2} }{2h^{2} } }
\end{array}\right].
$$

Impose the conditions
\begin{equation}
\label{eq_43}
h^{2} \le 4\tau \sigma ^{2} \,\,\, {\rm and} \,\tau \, \vert{\alpha^i_{k,j}}\vert 
\le {h/4}. \, \, \, \,\,
\end{equation}

\textbf{Remark 1}. Condition \eqref{eq_43} guarantees that all $\gamma^i_{k,j+1/2}$ are positive, and matrices ${\mathbb A}$ and ${\mathfrak A^i}$ are M-matrices \cite{Ple}.

\textbf{Proposition 1.} \textit{Components $m^i_{k,j+1/2}$ are non-negative $\forall i\in \{S,I,R,C\},\; \forall k = 1,...,M,\; \forall j=1,..., N-1 $ when $m^i_0(x_{j+1/2})\geq 0$ $\forall j=1,..., N-1$ and \eqref{eq_43} is performed.} 

\textit{\textbf{Proof.}} For non-negative components of $f^i_{k,\cdot}$ the non-negativity of $m^i_{k,\cdot}$ is provided by the monotonously property of the M-matrix \cite{Ple}. Put that $\exists\, i_0,\,k_0,\, j_0 $ for which $f^{i_0}_{{k_0-1},{j_0}}$ is negative and $m^{i_0}_{{k_0},{j_0}}$ change its value to negative one. From non-negativity of $m^{i_0}_0(x_{j+1/2}) \, \forall j=1,..., N-1$ follows that $\exists \zeta \in [(k_0-1)\tau,k_0\tau]$ that $m^{i_0,h}(\zeta,x_{j_0+1/2}) =0$. Then  $m^{i_0}_{{k_0-1},{j_0}} = O(\zeta - (k_0-1)\tau)$ and from \eqref{eq_41} it follows that $f^{i_0}_{{k_0-1},{j_0}}$ is non-negative for non-negative parameters $\beta,\,\mu,\,\varepsilon,\,\delta$. The obtained contradiction proves the proposition. 

\textbf{Remark 2}. Sum expressions in \eqref{eq_39} over $j= 0, ..., N-1 $ and multiply obtained expression by $\tau h$. As a result, for non-negative $m^{i,h}(t_{k-1},x)$ we get the equality  
\begin{equation}
\label{eq_44}
\begin{aligned}
\int_0^1 &\sum_{i=S,I,R,C} m^{i,h}(t_k,x) dx = \int_0^1 m^{h}(t_k,x) dx  =\\
&= \int_0^1 \sum_{i=S,I,R,C} m^{i,h}(t_{k-1},x) dx = \int_0^1 m^{h}(t_{k-1},x) dx. 
\end{aligned}
\end{equation}
since the following property is performed $\forall i \in \{S,I,R,C\}$
\begin{equation}
\label{eq_Assesm1_1}
 \gamma _{k,j+{1}/{2}\;}^{i,1} + 
 \gamma _{k,j+{1}/{2}\;}^{i,2}  + 
 \gamma _{k,j+{1}/{2}\;}^{i,3} = 1/\tau.   
\end{equation}
 Here $m^h$ is grid value of total mass of population and expression \eqref{eq_44} for non-negative $m^{i,h}$ is the discrete analogue of conservation law for total mass of agents.    
 
 \textbf{Proposition 2.} \textit{For \eqref{eq_39}--\eqref{eq_41} with initial \eqref{eq_FPK_init} and boundary \eqref{eq_30} conditions the following assessment is performed }
 \[\underset{0\le k\le M}{\mathop{\max }}\,{{\left\| {{m}^{i,h}}\left( {{t}_{k}},\cdot  \right) \right\|}_{1,h}}\le {{\left\| {{m}^i_{0}}\left( \cdot  \right) \right\|}_{1,h}}+T\underset{0\le k\le M}{\mathop{\max }}\,{{\left\| {{f }^{i,h}}\left( {{t}_{k}},\cdot  \right) \right\|}_{1,h}},\] 
 where ${{\|{{m}^{i,h}}\left( {{t}_{k}},\cdot  \right) \|}_{1,h}}$ is discrete analogue of $L_{1}(0,1)$–norm for grid function 
\[\int\limits_{0}^{1}{{{m}^{h}}\left( {{t}_{k}},x \right)\text{d}x}={{\left\| {{m}^{h}}\left( {{t}_{k}},\cdot  \right) \right\|}_{1,h}}:=\sum\limits_{j=0}^{N-1}{\vert m_{k,j+{1}/{2}}^{i,h} \vert h}.
\]

\textit{\textbf{Proof.}}
Use the key property \eqref{eq_Assesm1_1} of the approximation scheme \eqref{eq_39}. Put that $\Tilde{m}^{i,h}(t_k,\cdot)$ is the solution of \eqref{eq_39} for all non-negative $f^{i,h}_{k-1,j+1/2}$. Components $\Tilde{m}^{i,h}_{k,j+1/2}$ will be non-negative due to M-matrices properties \cite{Ple}. Multiply \eqref{eq_39} by $\tau$ and $h$ and sum up over $j=0,...,N-1$ for non-negative $f^{i,h}_{k-1,j+1/2}$. Taking into account \eqref{eq_Assesm1_1} we get an equality
 \[{{\left\| {\Tilde{m}^{i,h}}\left( {{t}_{k}},\cdot  \right) \right\|}_{1,h}} =  {{\left\| {\Tilde{m}^{i,h}}\left( t_k, \cdot  \right) \right\|}_{1,h}}+\tau{{\left\| {{f }^{i,h}}\left( {{t}_{k}},\cdot  \right) \right\|}_{1,h}}\]
 $\forall i\in\{S,I,R,C\},\;\forall k = 1,...,M$. Using mathematical induction on $k$ leads to the equality
 \[{{\left\| {\Tilde{m}^{i,h}}\left( {{t}_{k}},\cdot  \right) \right\|}_{1,h}} =  {{\left\| {\Tilde{m}^i_{0}}\left( \cdot  \right) \right\|}_{1,h}}+\tau\sum_{s=1}^{k}{{\left\| {{f }^{i,h}}\left( {{t}_{s}},\cdot  \right) \right\|}_{1,h}}.\]

Now put $\Tilde{m}^{i,h}(t_k,\cdot)$ the solution of system \eqref{eq_39}, where $f^{i,h}_{k-1,j+1/2}$ can be of any sign. Note that 
\[ \Tilde{m}^{i,h}\left( {{t}_{k}},x_{j+1/2} \right) - {m}^{i,h}\left( {t}_{{k},{j+1/2}}  \right)  \geq 0. \]
It leads from monotone property of M-matrix \cite{Ple} after substitution the expression $\Tilde{m}^{i,h}\left( {{t}_{k}},x_{j+1/2} \right) - {m}^{i,h}\left( {}{{t}_{k}},_{j+1/2}  \right) $ to \eqref{eq_39}. From non-negativity of components $\Tilde{m}^{i,h}_{k,j+1/2}$ and ${m}^{i,h}_{k,j+1/2}$ it follows that 
 \begin{equation}
 \label{eq_Assesm1_2}
 {\left\| {{m}^{i,h}}\left( {{t}_{k}},\cdot  \right) \right\|}_{1,h} \leq {{\left\| {\Tilde{m}^{i,h}}\left( {{t}_{k}},\cdot  \right) \right\|}_{1,h}} =  {{\left\| {\Tilde{m}^i_{0}}\left( \cdot  \right) \right\|}_{1,h}}+\tau\sum_{s=1}^{k}{{\left\| {{f }^{i,h}}\left( {{t}_{s}},\cdot  \right) \right\|}_{1,h}}.   
 \end{equation}
Taking the maximum of both parts of \eqref{eq_Assesm1_2}, we obtain the required estimate. 

\subsection{Discrete optimal control problem }
 
 Instead of integral cost function in \eqref{eq_11} use the discrete one
 \begin{equation}
\label{eq_47}
\begin{aligned}
   & {{J}^{h}}({{m}^{h}_{SIRC}},{{{\alpha }}^{h}_{SIRC}})= \sum\limits_{i}\sum\limits_{k=0}^{M-1}{\sum\limits_{j=0}^{N-1}{\left. \left( r_{k,j+{1}/{2}\;}^{i,h}m_{k,j+{1}/{2}\;}^{i,h} \right.+g_{k,j+{1}/{2}\;}^{i,h} \right)\tau }}h+\\ &+\sum\limits_{j=0}^{N-1}h{(m_{M,j+{1}/{2}\;}^{I,h})^2/2}.
\end{aligned}
\end{equation}
Here $r_{k,i+{1}/{2}\;}^{h}$ is carried out for $F\left( {\alpha },t,x \right)$ by the following way:
\begin{equation}
\label{eq_48}
r_{k,j+{1}/{2}\;}^{i,h}={F^i\left({\alpha }_{k,j}^{h},{{t}_{k}},{{x}_{j}} \right)}/{2}\;+\,{F^i\left( {\alpha }_{k,j+1}^{h},{{t}_{k}},{{x}_{j+1}} \right)}/{2}\; 	    
\end{equation}
with  $\alpha _{k,j}^{i,h}:={{\alpha }^{i,h}}\left( {{t}_{k}},{{x}_{j}} \right).$ Then instead of differential minimization problem we get the grid one:
\begin{equation}
\label{eq_49}
\left\{ \begin{aligned}
  & \underset{\alpha_{SIRC} }{\mathop{\inf }}\,{{J}^{h}}({{m}^{h}_{SIRC}},{{{\alpha }}^{h}_{SIRC}}), \\ 
 & \mathfrak{A^i}m_{\cdot ,\cdot }^{i,h}=\mathfrak{F^i}m_{0,\cdot }^{i,h}. \\ 
\end{aligned} \right.
\end{equation}

 \textbf{Remark 3.} The decomposition in a series of Taylor shows that the finite difference problem \eqref{eq_49} approximates the initial differential formulation  with the order of $O(\tau+h^2).$

Introduce also the following form
\begin{equation}
\label{eq_50}
\left\langle a,b \right\rangle =h\sum\limits_{j=0}^{N-1}{{{a}_{j+{1}/{2}\;}}{{b}_{j+{1}/{2}\;}}}
\end{equation} 
for any vectors $a=\left\{ a_{j+1/2}\right\}_{j=0,...,N-1}$ and $b=\left\{ b_{j+1/2}\right\}_{j=0,...,N-1}$.
To formulate a discrete optimal control problem, introduce into consideration the set of grid functions $ v _ {\cdot, \cdot} ^ {i,h} = \left\{v_{k, j + {\text{1}}/{\text{2}}\;}^{i,h}\right\}_{j = 0, ..., N-1}^{k = 0, ..., M}$. Multiply the $k-$th component of $i-$th equation from \eqref{eq_42} by $ \tau v_ {k,\cdot}^ {i,h}$ and sum over $k$ taking into account the notation \eqref{eq_50}:
\begin{equation}
\label{eq_51}
\begin{aligned}
  & L^{i,h}: =  \tau \sum\limits_{k=1}^{M}{\left( \left\langle \mathbb{A}m_{k,\cdot }^{i,h},v_{k,\cdot }^{i,h} \right\rangle -\left\langle {{\mathbb{B}}^{i}_{k}}m_{k-1,\cdot }^{i,h},v_{k,\cdot }^{i,h} \right\rangle   - \langle v_{k,\cdot }^{i,h},f_{k-1,\cdot }^{i,h} \rangle \right)} +\\
  & =-\left\langle m_{0,\cdot }^{i,h},\mathbb{A}v_{0,\cdot }^{i,h} \right\rangle \tau  +\left\langle m_{M,\cdot }^{i,h},\mathbb{A}v_{M,\cdot }^{i,h} \right\rangle \tau +\\
  &+\tau \sum\limits_{k=0}^{M-1}{\left( \left\langle m_{k,\cdot }^{i,h},\mathbb{A}v_{k,\cdot }^{i,h} \right\rangle -\left\langle m_{k,\cdot }^{i,h},(\mathbb{B}_{k+1}^{i})^{*}v_{k+1,\cdot }^{i,h} \right\rangle - \langle v_{k+1,\cdot }^{i,h},f_{k,\cdot }^{i,h} \rangle \right)}. 
\end{aligned}
\end{equation}
Now write down the Lagrangian for the grid optimization problem \eqref{eq_49}
\begin{equation}
\label{eq_52}
\begin{aligned}
  & {{\Im }^{h}}({{m}_{SIRC}^{h}},{{\alpha }_{SIRC}^{h}},{{v}_{SIRC}^{h}}):={{J}^{h}}({m_{SIRC}^{h}},{{\alpha }_{SIRC}^{h}})-\sum_{i = \{S,I,R,C\}} L^{i,h}. \\ 
\end{aligned}
\end{equation}
Here $ (\mathbb {B}_{k}^{i})^{*} = (\mathbb{B}_{k}^{i})^{T} $ means the matrix conjugate to $ \mathbb{B}_{k}^{{i}} $. Then the problem of finding a saddle point \eqref{eq_22} in the grid case takes the form
\begin{equation}
\label{eq_53}
\underset{({{m}_{SIRC}^{h}},{{\alpha }_{SIRC}^{h}})}{\mathop{\inf }}\,\,\,\underset{{{v}_{SIRC}^{h}}}{\mathop{\sup}}\,\,\,{{\Im }^{h}}({{m}_{SIRC}^{h}},{{\alpha }_{SIRC}^{h}},{{v}_{SIRC}^{h}}).
\end{equation}

Differentiate Lagrangian (\ref{eq_52}) with respect to the individual components. We get the following system of algebraic equation 
\begin{equation}
\label{eq_54}
\mathbb{A}v_{k,\cdot }^{i,h}=(\mathbb{B}_{k+1}^{i})^{*}v_{k+1,\cdot }^{i,h}+z_{k,\cdot }^{i,h} 
\end{equation}
$\forall \,i\in \{S,I,R,C\}, \forall \,k=M-1,M-2,...,0.$ Here
\begin{equation}
\label{eq_55}
\begin{aligned}
z_{k,j+{1}/{2}\;}^{i,h}=\partial g^i /\partial m^i &\left( {{t}_{k}},{{x}_{j+{1}/{2}\;}},m_{k,j+{1}/{2}\;}^{h} \right)+r_{k,j+{1}/{2}\;}^{i,h}  \\
&+ \sum_{l=\{S,I,R,C\}}{v^{l,h}_{k+1}\partial f^{l,h}_{k} / \partial m^{i,h}_{k,\cdot}}
\end{aligned}
\end{equation}
and
\begin{equation}
\label{eq_56}
\hspace*{-43mm}v_{k,-{1}/{2}\;}^{h}=v_{k,{1}/{2}\;}^{h},\text{ }v_{k,N+{1}/{2}\;}^{h}=v_{k,N-{1}/{2}\;}^{h}
\end{equation}
$\forall \,\,k=M-1,...,0\text{  }\forall \,\,i=0,...,N-1.$
The initial conditions for \eqref{eq_54} after variation of the Lagrangian \eqref{eq_52} can be written in the following form
\begin{equation}
\label{eq_57} 
\left\{
\begin{aligned}
&v_{M,i+{1}/{2}\;}^{i,h}=0,\;\; \text{if } i\in \{S,R,C\},\\
&\mathbb{A}v_{M,i+{1}/{2}\;}^{I,h}=m^{I,h}_{M,i+{1}/{2}\;} \; \text{if } i=I.  
\end{aligned}
\right.
\end{equation}
\textbf{Proposition 3.} \textit{For \eqref{eq_54}--\eqref{eq_57} under the restrictions \eqref{eq_43} the following assessments are performed }
 \[\underset{0\le k\le M}{\mathop{\max }}\,{{\left\| {{v}^{i,h}}\left( {{t}_{k}},\cdot  \right) \right\|}_{\infty,h}}\le T\underset{0\le k\le M}{\mathop{\max }}\,{{\left\| {{z }^{i,h}}\left( {{t}_{k}},\cdot  \right) \right\|}_{\infty,h}}\] 
 for $\in \{S,R,C\}$ and for $i=I$:
  \[\underset{0\le k\le M}{\mathop{\max }}\,{{\left\| {{v}^{I,h}}\left( {{t}_{k}},\cdot  \right) \right\|}_{\infty,h}}\le {{\left\| {{m}^I(t_M,\cdot)} \right\|}_{\infty,h}}+T\underset{0\le k\le M}{\mathop{\max }}\,{{\left\| {{z }^{i,h}}\left( {{t}_{k}},\cdot  \right) \right\|}_{\infty,h}},\] 
 where ${{\|{{m}^{i,h}}\left( {{t}_{k}},\cdot  \right) \|}_{\infty,h}}$ is discrete analogue of $L_{\infty}(0,1)$–norm for grid function 
\[\underset{0\le x\le 1}{\mathop{\max }}\,\vert {{v}^{h}}\left( {{t}_{k}},x \right) \vert={{\left\| {{v}^{h}}\left( {{t}_{k}},\cdot  \right) \right\|}_{\infty ,h}}:=\underset{0\le i\le N-1}{\mathop{\max }}\,\vert v_{k,i+{1}/{2}\;}^{h} \vert.\]
\textit{\textbf{Proof.}} Put $\vert \Tilde v^{i,h}(t_k, x_{j+1/2}) \vert$ as component which is  maximal in absolute value on the layer $t_k$ so that $\vert \Tilde v^{i,h}(t_k, x_{j+1/2}) \vert =\| v^{h,i}(t_k,\cdot) \|_{\infty,h}$. Use again  the key property of the coefficients \eqref{eq_Assesm1_1} to obtain the inequality
\[\| v^{h,i}(t_k,\cdot) \|_{\infty,h} = \vert \Tilde v^{i,h}(t_k, x_{j+1/2}) \vert
\leq  \| v^{h,i}(t_{k+1},\cdot) \|_{\infty,h} +\ \tau \| z^{h,i}(t_k,\cdot) \|_{\infty,h}.  \]
The use of mathematical induction on $k$ leads to the estimate
\begin{equation}
\label{eq_Asessm2_1}
\| v^{h,i}(t_k,\cdot) \|_{\infty,h} 
\leq  \| v^{h,i}(t_{M},\cdot) \|_{\infty,h} +\ (M-k)\tau \| z^{h,i}(t_k,\cdot) \|_{\infty,h}.   
\end{equation}
For $i \in \{S,R,C\}$ with zero "initial" conditions for \eqref{eq_54} the required estimate can be obtained after taking a maximum over $k$ in \eqref{eq_Asessm2_1}. For $i=I$ from \eqref{eq_57} follows that 
$$ \| v^{h,i}(t_{M},\cdot) \|_{\infty,h} \leq  \| \mathbb{A}^{-1}\|_{\infty,h} \| m^{h,I}(t_M,\cdot) \|_{\infty,h}.
$$
It is known that for an M-matrix with strict diagonal dominance the following statement holds (Theorem 2 from \cite{Vol}): 
$$ \| \mathbb{A}^{-1}\|_{\infty,h}  = 1/R,
$$
where $R$ is the amount of diagonal dominance, which is the same for each row of the matrix $\mathbb{A}$ and equaled to $1/\tau$. Then for $i=I$ \eqref{eq_Asessm2_1} can be rewritten as 
\begin{equation}
\label{eq_Asessm2_2}
\| v^{h,I}(t_k,\cdot) \|_{\infty,h} 
\leq  \tau \| m^{h,I}(t_{M},\cdot) \|_{\infty,h} +\ (M-k)\tau \| z^{h,I}(t_k,\cdot) \|_{\infty,h}.
\end{equation}
Taking a maximum over $k$ in \eqref{eq_Asessm2_2} we get the required estimate for $i=I$.

Also \eqref{eq_53} gives the following grid analogue for optimality conditions \eqref{eq_29} in point $\left( {{t}_{k}},{{x}_{j}} \right):$ 
\begin{equation}
\label{eq_58}
\frac{\partial F_i}{\partial \alpha_i }\left( \alpha _{k,j}^{i,h},{{t}_{k}},{{x}_{j}} \right)+\frac{v_{k,j+{1}/{2}\;}^{i,h}-v_{k,j-{1}/{2}\;}^{i,h}}{h}=0
\end{equation} 
$\,\forall \,k=1,\ldots ,M\;\forall \,j=1,\ldots ,N-1$.

Thus, the solution to the discrete optimization problem \eqref{eq_49} can be found by the following \textbf{iterative algorithm}: 
\begin{enumerate}
    \item Put the initial value of grid functions $\alpha^{i,h}_{\cdot,\cdot} $ equal to zero.
    \item For zero functions $\alpha^{i,h}_{\cdot,\cdot} $ get the initial value of functions $m^{i,h}_{\cdot,\cdot}$ using \eqref{eq_39}-\eqref{eq_42} and $J^{h}(m^{i,h}_{\cdot,\cdot},\alpha^{i,h}_{\cdot,\cdot})$ using \eqref{eq_47}, \eqref{eq_48}.
    \item Get the value of $v^{i,h}_{\cdot,\cdot}$ functions solving the system \eqref{eq_54}-\eqref{eq_57}.
    \item Get a new value for functions $\alpha^{i,h}_{\cdot,\cdot}$ by the solution of \eqref{eq_58}.
    \item Get a new value for functions $m^{i,h}_{\cdot,\cdot}$ by the solution of \eqref{eq_39}--\eqref{eq_43}.
    \item Get a new value cost function $J^{h}(m^{i,h}_{\cdot,\cdot},\alpha^{i,h}_{\cdot,\cdot})$ by \eqref{eq_47}, \eqref{eq_48}.
    \item If the cost function $J^{h}(m^{i,h}_{\cdot,\cdot},\alpha^{i,h}_{\cdot,\cdot}) $ reaches its minima with the given accuracy then choose obtained on the previous steps functions $m^{i,h}_{\cdot,\cdot},\alpha^{i,h}_{\cdot,\cdot}$ as a solution of \eqref{eq_49}. Otherwise go to step 3 for new iteration. 
\end{enumerate}

\subsection{Discrete MFG with corrective control}
For discrete statement of optimal control problem with an external influence on agent's strategy we can use schemes \eqref{eq_39}-\eqref{eq_43} and \eqref{eq_54}-\eqref{eq_57} introducing into consideration the following representation  
\begin{equation}
\label{eq_59}
\alpha^{i,h}_{k,j} = \hat{\alpha }^{i,h}_{k,j}+\rho^{h} \left( \hat{\alpha }^{i,h}_{k,j} \right)
\end{equation}
$\forall i\in\{S,I,R,C\},\;\forall k=1,...,M\;\forall j = 1,...,N-1$. Instead of integral function \eqref{eq_28} use the discrete one 
 \begin{equation}
\label{eq_60}
\begin{aligned}
   & {{J}^{h}}({{m}^{h}_{SIRC}},{{\hat{\alpha }}^{h}_{SIRC}})= \sum\limits_{i}\sum\limits_{k=0}^{M-1}{\sum\limits_{j=0}^{N-1}{\left. \left( \hat{r}_{k,j+{1}/{2}\;}^{i,h}m_{k,j+{1}/{2}\;}^{i,h} \right.+g_{k,j+{1}/{2}\;}^{i,h} \right)\tau }}h+\\ &+\sum\limits_{j=0}^{N-1}h{(m_{M,j+{1}/{2}\;}^{I,h})^2/2},
\end{aligned}
\end{equation}
where $\hat{r}_{k,i+{1}/{2}\;}^{h}$ is carried out for $F\left( \hat{\alpha },t,x \right)$ by the following way:
\begin{equation}
\label{eq_61}
\hat{r}_{k,j+{1}/{2}\;}^{i,h}={F^i\left({\hat{\alpha} }_{k,j}^{h},{{t}_{k}},{{x}_{j}} \right)}/{2}\;+\,{F^i\left( \hat{\alpha }_{k,j+1}^{h},{{t}_{k}},{{x}_{j+1}} \right)}/{2}.\; 	    
\end{equation}
Then we can write the grid optimization problem in the following form:
\begin{equation}
\label{eq_62}
\left\{ \begin{aligned}
  & \underset{\hat{\alpha}_{SIRC} }{\mathop{\inf }}\,{{J}^{h}}({{m}^{h}_{SIRC}},{{{\hat{\alpha} }}^{h}_{SIRC}}), \\ 
 & \mathfrak{A^i}m_{\cdot ,\cdot }^{i,h}=\mathfrak{F^i}m_{0,\cdot }^{i,h}. \\ 
\end{aligned} \right.
\end{equation} 
The solution of \eqref{eq_62} can be found by iterative algorithm described in previous section taking into account the representation \eqref{eq_59}. 

\subsection{Initial dataset} \label{initial_data}
For the numerical implementation the constructed algorithm we use the following parameters. In table \ref{tabEpid} the epidemiological constants are described. The presented values are obtained by the solution the inverse problem for the corresponding statistical data on the incidence of COVID-19 in the Novosibirsk region, Russia \cite{Cha}.  
\begin{table}[h]
  \caption{Epidemiological constants description and its reference values for Novosibirsk region}
  \centering
 \begin{tabular}{|p{135pt}|p{30pt}|p{40pt}|p{40pt}|p{40pt}|}
\hline

\centering \textbf{Parameter description} & \textbf{Symbol} & \textbf{05/01/20 -- -- --  06/30/20}  & \textbf{06/30/20 -- -- -- 08/08/20} & \textbf{12/01/20 -- -- -- 03/10/21}\\

\hline
The contact/infection transmission rate & $\beta$ & 0.2821 & 0.3253 & 0.4145\\
 \hline
 The recovery rate of the infected population & $\gamma$ & 0.2530 & 0.3466 &  0.4257 \\
 \hline
The rate at which the recovered population becomes the cross-immune population and moves from total to partial immunity & $\delta$ & 0.3657 & 0.0794 & 0.0889\\
 \hline
The rate at which the cross-immune population becomes susceptible again & $\mu$ & 0.0060 &  0.0376& 0.0267\\
 \hline
The average reinfection probability of a cross-immune individual & $\varepsilon$ & 0.0361 & 0.1446 & 0.0928\\
 \hline
\end{tabular}
  \label{tabEpid}
\end{table}

In addition for MFG model we use the parameters presented in table \ref{tabMFGconst}.
\begin{table}[h]
  \caption{Description of constants of MFG model }
  \centering
 \begin{tabular}{|p{195pt}|c|c|}
\hline
\centering\textbf{Parameter description} & \textbf{Symbol} & \textbf{Value} \\
 \hline
 Time horizon & $T$ & 100 days \\
 \hline
 The stochastic parameter of the system. Diffusion in physical meaning & $\sigma^2_i$ & 0.2 \\
 \hline
 Coefficient characterizing the degree of growth in the costs of an individual for non-compliance with physical distancing  &$ c_{1i}$  &  
$ \left\{\begin{array}{l}
6 \; \text{ if } i \in \{S,R,C\}, \vspace{2mm}\\
2 \;\text{ if } i =I
\end{array}\right.
$  \\
 \hline
 The coefficient of influence of the entire mass of agents on the opinion of an individual, $c_{2i}\in [0,1]$ & $c_{2i}$ &  0.9\\
  \hline
Penalty for deviating the agent from the social distancing strategy, $c_{3i}\in [0,2]$ & $c_{3i}$ &
$ 0.7 $ \\
 \hline
\end{tabular}
  \label{tabMFGconst}
\end{table}
Note that combination of constants $c_{1i}$ and $c_{2i}$ is quite close in physical meaning to the so-called "index of self-isolation" which was introduced by Yandex company (https://yandex.ru/company/researches/2020/podomam). The current values of constants are corresponding to the low level of self-isolation (in terms of Yandex it is about 0.5 -- 1), that corresponds the situation in time periods under consideration. The values of the coefficients $c_{3i}$ are chosen in this way because  in these time periods the quarantine measures introduced by the government were not so restrictive.  

And finally as an initial distributions of mass of agents we put 
\begin{equation}
    \label{eq_init}
m_{0i} = \frac{A_i}{B_i} \left( exp\left(-\frac{(x-x_i^c)^2}{2 (\sigma^c_i)^2}\right)\bigg/{\sigma^c_i\sqrt{2\pi}}  + a_ix^2 + b_i(1-x)^2 \right),
\end{equation}
where $A_i$ is proportion of current fraction in relative to total population at initial time moment; $B_i$ is normalization factor equaled to integral over $[0,1]$ from expression in brackets; $a_i = exp\left(-\frac{(1-x_i^c)^2}{2 (\sigma^c_i)^2}\right) (1-x^c_i)/(2(\sigma^c_i)^3\sqrt{2\pi})$ and $b_i = exp\left(-\frac{(x_i^c)^2}{2 (\sigma^c_i)^2}\right) (x^c_i)/(2(\sigma^c_i)^3\sqrt{2\pi})$ to ensure boundary conditions \eqref{eq_7} for $m_{0i}$. We use the following values for  $x^c_i$ and $\sigma^c_i$
\begin{equation}
\label{eq_65}
\begin{aligned}
   & x_S^c = 0.8 ;\;\;\sigma_S^c = 0.1; \\
   & x_I^c = 0.2;\;\;\sigma_I^c =0.1 ; \\
   & x_R^c = 0.7;\;\;\sigma_R^c =0.2 ;\\
   & x_C^c =0.3 ;\;\;\sigma_C^c = 0.2 .\\
\end{aligned}
\end{equation}
The physical meaning of constants \eqref{eq_65} is that non-infected (S) people do not seek to comply with restrictions and self-isolation for economic reasons; those infected (I) in general comply with quarantine measures. The recovered part of the population (R) has full immunity and often does not comply with social measures, since they are confident in their immunity; and cross-immune ones (C) are afraid of re-infection after a while.

For $A_i$ constants we choose the following parameters, presented in table \ref{tabInit} for different time periods.
\begin{table}[!ht]
  \caption{Initial data for population distribution}
  \centering
 \begin{tabular}{|c|c|c|c|}
\hline
\multirow{2}*{\textbf{Symbol}} & \multicolumn{3}{c|}{\textbf{Time period}}\\
\cline{2-4}
 & \centering\textbf{05/01/20 -- 06/30/20}  & \centering\textbf{06/30/20 -- 08/08/20} & \textbf{12/01/20 -- 03/10/21}\\
 \hline
 $A_S$ & 0.999757 &  \multirow{4}{80pt}{The data obtained in horizon time $T$ from previous time period} & 0.991199\\
 $A_I$  & 0.000204 &   & 0.001505 \\
 $A_R$ & 0.000039 &  & 0.006937\\
$A_C$ & 0.000000 &  & 0.000359\\
 \hline
\end{tabular}
  \label{tabInit}
\end{table}

Here parameters $A_I$ are also the solution of the inverse problem. 

\subsection{Numerical results for modelling of COVID-19 propagation in Novosibirsk region, Russia}
\begin{figure}[!ht]
    \centering
    \includegraphics[width=0.9\textwidth]{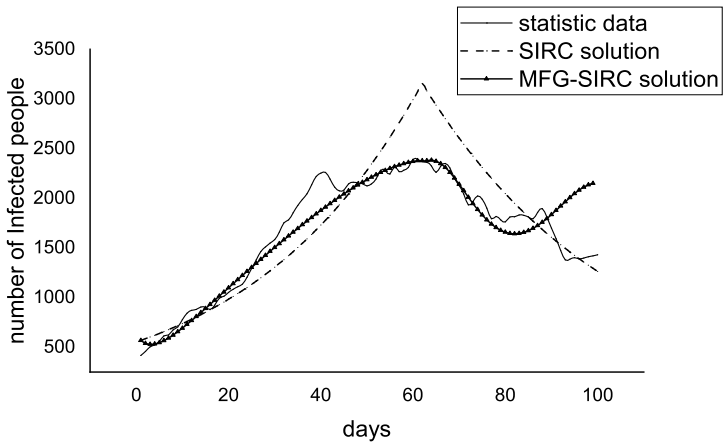}
    \caption{The comparison of SIRC and MFG-SIRC models for number of daily diagnoses people in Novosibirsk region.}
    \label{figMay}
\end{figure}

First, let us show the difference between using the standard SIRC differential model and the MFG model based on it. The comparison will be made with the official statistical data of the coronavirus report \cite{Cha} for a period of 100 days from May 1, 2020 in Novosibirsk city, Russia. The collected statistic data, as well the solutions of differential SIRC and grid MFG-SIRC models are presented in figure \ref{figMay}. 

As it is clearly seen, the statistical data for the period under consideration have two areas: the rise and fall of the incidence. Since models of the SIR type are poorly applied for long-term prediction, here, to approximate the epidemiological situation, two curves with different epidemiological parameters were glued together. The epidemiological parameters of the model were restored according to statistical data and the optimization method presented in the work~\cite{Krivorotko_2020} for two time intervals: from May 1 to June 30 and from June 30 to August 8, 2020. The obtained parameters are presented in table \ref{tabEpid}. These parameters were used for numerical implementation both for the SIRC model and MFG-SIRC one. The spatial and stochastic constants for MFG implementation are presented in section \ref{initial_data}. As can be seen from the figure, MFG gives a more accurate approximation to statistical data, but like the parent SIRC model, it does not work well for long-term forecasts.

\begin{figure}[!ht]
    \centering
    \includegraphics[width=0.9\textwidth]{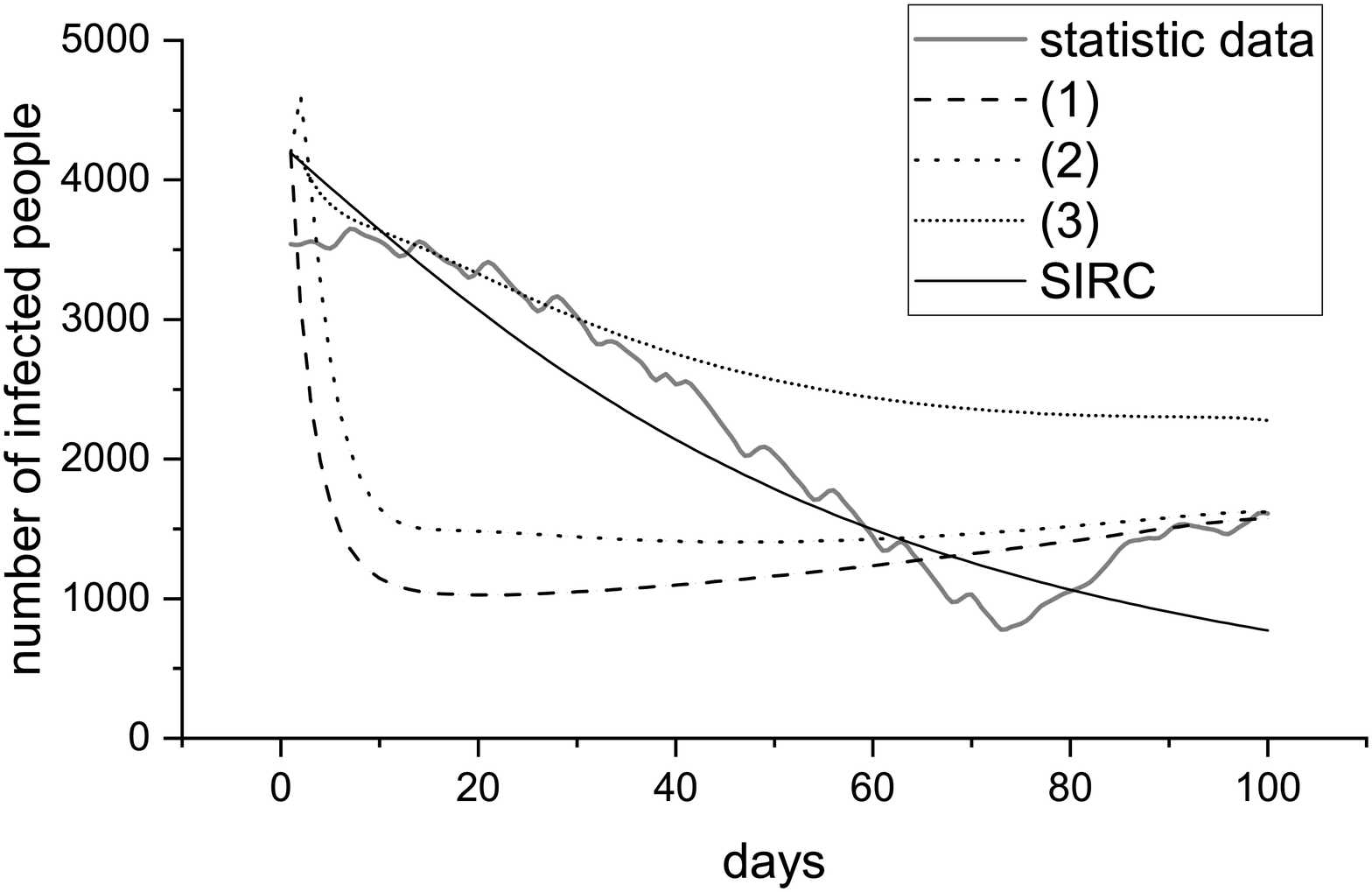}
    \caption{The comparison of MFG-SIRC with different initial distributions}
    \label{figDec}
\end{figure}

Figure \ref{figMay} also shows that introducing spatial distribution into consideration has a significant impact on population behavior. Figure \ref{figDec} shows the dependence of the obtained solution on the selected initial spatial distribution. The comparison was made on the basis of the coronavirus report in the Novosibirsk city in the period from December 1, 2020, when there was a decline in the incidence. For this time period, the restored epidemiological parameters are also presented in the table \ref{tabEpid}, the stochastic parameter of the system $\sigma^2$ was chosen equal 0.5, and the following curves were used as the initial distributions for different population groups:

(1) -- initial distributions $m_{0i}$ are determined by \eqref{eq_init} with parameters \eqref{eq_65};

(2) -- initial distributions $m_{0i}$ are determined by \eqref{eq_init} with parameters $ x_i^c = 0.2;\;\sigma_i^c =0.1 $ $ \forall i\in\{S,I,R,C\}$;

(3) -- initial distributions $m_{0i}$ are determined by \eqref{eq_init} with parameters $ x_i^c = 0.5;\;\sigma_i^c =0.2 $ $ \forall i\in\{S,I,R,C\}$. 

In physical sense, it means that at initial time in the second case (2) agents are clearly positive to quarantine measures and physical distancing, and in third case (3) are not so patient to performed it. The motivation for case (1) is presented in section \ref{initial_data}.

\begin{figure}[!ht]
    \centering
    \includegraphics[width=0.9\textwidth]{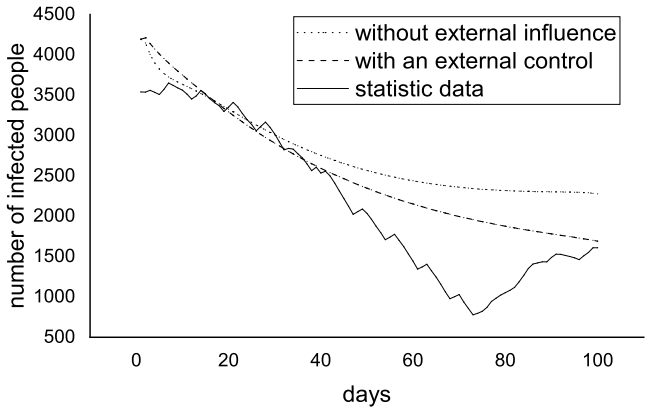}
    \caption{The comparison of MFG-SIRC with external influence and without it (case (3) in  figure~\ref{figDec}) }
    \label{figExt}
\end{figure}

And finally, for case (3), we introduce external corrective control as a solution of the system \eqref{eq_62}. Comparison of the obtained curves is shown in Figure \ref{figExt}.

The figure \ref{figExt} shows that taking into account external restrictive measures that are not the agent's choice has a significant impact on population dynamics.

\section{Conclusion and discussions}

The article is devoted to the transfer of well-known economic models of the "mean field" to forecasting the spread of epidemics, in particular, COVID-19. The use of the MFG  approach is due to the fact that traditional epidemiological models such as SIR and others based on it do not take into account the heterogeneity of the population, and therefore cannot be used for long-term forecasts. Another well-known approach to the spread of viruses, the so-called agent-based models, allow taking into account non-epidemiological factors, but lead to computationally complex systems \cite{Krivorotko_IDM_2021}. In turn, structural simplicity and taking into account the spatial characteristics of MFG models can solve both problems.

Numerical experiments on modeling the dynamics of the spread of COVID-19 in Novosibirsk, Russia have shown that despite the fact that the dynamics of the population is determined by the differential SIRC model, taking into account the spatial characteristics of agents has a huge impact on the final result. This gives an advantage over SIR-type models but also generates a huge class of problems for determining the parameters of the system, both from the field of statistical analysis and inverse problems. In addition, there are several more ways to modify the presented model. The greatest interest is the consideration of the epidemiological parameters $\beta,\gamma$ and others depending on the position $x$ of the agent in space or even on its strategy $\alpha(t,x)$. Note also that the proposed computational approach can be carried over to any model of the SIR type.

The numerical algorithm considered in this paper is based on the ideas proposed in the works \cite{Sha1, Sha2, Sha3, Sha4} and modified for application to the problems of epidemiology. The algorithm  gives a direct and simply rule for minimizing the cost functional, ensures the fulfillment of the law of conservation of the entire mass of agents, and allows to take into account more complex functions $F_i(\alpha,t,x)$ responsible for control, instead of the traditionally used quadratically dependent on $\alpha$ \cite{Ben}.

\section{Acknowledgement}
This work was supported by the Russian Science Foundation (project No.~18-71-10044).


\end{document}